\documentclass[reqno]{amsart}
\usepackage{amssymb}

\def\br#1{#1{}\\{}#1}
\def\R{\mathbb R}
\def\N{\mathbb N}
\def\P{\mathbb P}
\def\Z{\mathbb Z}
\def\T{\mathbf T}
\let\d=\delta
\let\g=\gamma
\let\G=\Gamma
\let\e=\varepsilon
\let\tu=\textup
\def\deg{\operatorname{deg}}

\newtheorem{thm}{Theorem} \newtheorem*{thm*}{Theorem}
\theoremstyle{definition}
\newtheorem*{dfn*}{Definition}
\newtheorem{pro}{Problem}
\newtheorem*{Turanpro}{Tur\'an's problem}
\newtheorem{exa}{Example}

\begin{document}

\title{Relation between Tur\'an extremum problem and~van~der Corput sets}
\author{D.V.~Gorbachev, A.S.~Manoshina (Tula)}
\email{dvg@uic.tula.ru; ann@mm.tula.ru}
\urladdr{http://home.uic.tula.ru/$\sim$gd030473}
\thanks{This research was supported by the RFBR under grant no.\;03-01-00647 and 03-01-06199.
\endgraf
Matem. zametki (Russia), 2004 (in publishing).}
\date{12/16/2003}
\begin{abstract}
Let $K\subset\mathbb N$ and $\mathbf T(K)$ is a set of trigonometric
polynomials
\[
T(x)=T_0+\sum_{k\in K,\,k\le H}T_k\cos(2\pi kx), \qquad H>1,
\]
$T(x)\ge0$ for all $x$ and $T(0)=1$.

Suppose that $0<h\le1/2$ and $K(h)$ is the class of functions
\[
f(x)=\sum_{n=0}^{\infty}a_n\cos(2\pi nx)
\]
satisfying the following conditions: $a_n\ge0$ for all $n$, $f(0)=1$ and $f(x)=0$ for
$h\le|x|\le1/2$.

We consider an relation between extremum problem
\[
\delta(K)=\inf_{T\in\mathbf T(K)}T_0
\]
and Tur\'an extremum problem
\[
A(h)=\sup_{f\in K(h)}a_0=\sup_{f\in K(h)}\int_{-h}^hf(x)\,dx
\]
for rational numbers $h=p/q$ and set $K=\bigcup\limits_{\nu=0}^\infty\{q\nu+p,\ldots,q\nu+q-p\}$.

The problem $\delta(K)$ is connection with van der Korput sets. Van der Korput sets study
in analytic number theory.
\end{abstract}
\maketitle

In number theory following question is important (see \cite{Mon}, we use results of
Chapter~2): is given set of numbers
\[
\{u_n\}_{n=1}^\infty\subset\R
\]
uniformly distributed?

To date several methods are known for solving this problem \cite{Mon}. One way was
given by J.G.~van der Corput. This method consists in research of sequence
\[
\{u_{n+k}-u_n\}_{n=1}^\infty,
\]
where $k$ runs over some set $K$ of natural numbers.

\begin{dfn*}[\cite{Mon}]
A set $K\subset\N$ is called a van der Corput set if the sequence
$\{u_n\}_{n=1}^\infty$ is  uniformly distributed (mod\;1) whenever the differenced
sequence
\[
\{u_{n+k}-u_n\}_{n=1}^\infty
\]
is uniformly distributed (mod 1) for all $k\in K$.
\end{dfn*}

Set of natural numbers $K=\N$ is a van der Corput set. It was proved by van der
Corput \cite{Mon}.

In [1], criterion was given for determining whether a set $K\subset\N$ is a van der
Corput set. It is founded on following extremum problem for positive trigonometric
polynomials.

\begin{pro}\label{td}
Let $\T(K)$ be a set of trigonometric polynomials $T(x)$ such as following
conditions are satisfied:

\smallskip
1) $T(x)=T_0+\sum\limits_{\substack{k\in K\\ k\le H}}T_k\cos(2\pi kx)$, $H>1$;

\smallskip
2) $T(x)\ge0$ $\forall x\in\R$;

\smallskip
3) $T(0)=1$.

\smallskip
One needs to find value
\[
\d(K)=\inf_{T\in\T(K)}T_0.
\]
\end{pro}

\begin{thm*}[\cite{Mon}]
A set $K\subset\N$ is a van der Corput set if and only if
\[
\d(K)=0.
\]
\end{thm*}

From this it follows that $\d(\N)=0$ (van der Cortup's result) and $\d(K)>0$ for
any finite set $K$. Therefore a van der Corput set is not finite.

Exact value of $\d(K)$ is known in few cases. Now, two examples of exact
values~\cite{Mon}.

\begin{exa}\label{p1}
Suppose $q\in\N$, $q\ge2$,
\[
K_q^0=\{1,2,\ldots,q-1\},\qquad
K_q=q\Z_++K_q^0=\{k\in\N:q\nmid k\},
\]
here $q\Z_++K=\{q\nu+k:\nu\in\Z_+,\ k\in K\}$. Then
\[
\d(K_q^0)=\d(K_q)=\frac1q.
\]
Thus this set $K_q$ is not a van der Corput set.
\end{exa}

\begin{exa}\label{p2}
$\d(\{2,3\})=\dfrac{\cos(\pi/5)}{1+\cos(\pi/5)}=0{,}44721\ldots$
\end{exa}

Many properties are established for value $\d(K)$ in \cite{Mon}. Now, some of them
(here $K,K_1,K_2\subset\N$, $q\in\N$):

\smallskip
1) if $K_1\subset K_2$, then $\d(K_1)\ge \d(K_2)$;

\smallskip
2) $\d(qK)=\d(K)$, where $qK=\{qk:k\in K\}$;

\smallskip
3) $\d(K^{(q)})\le q\d(K)$, where $K^{(q)}=\{k\in K:q\mid k\}$ (if set $K^{(q)}$ is
empty, then $\d(K)\ge1/q$);

\smallskip
4) $\d(K_1)\d(K_2)\le\d(K_1\cup K_2)$.

\smallskip
These properties show that it is desirable to know value $\d(K)$ for given set $K$
even if ${\d(K)\ne 0}$. Examples from \cite{Mon}: set $Q=\{\nu^2+1\}_{\nu=1}^\infty$
and set of prime numbers~$\P$ are not van der Corput sets, since $Q\subset K_3$,
$\P\subset K_4$ and $\d(Q)\ge1/3$, $\d(\P)\ge1/4$.

There exists a relation between extremum problem $\d(K)$ and Tur\'an extremum problem
$A(h)$ \cite{GorMan} for rational numbers $h=p/q$.

\begin{Turanpro}[\cite{GorMan}]
Suppose that $0<h\le1/2$ and $K(h)$ is the class of continuous 1-periodic even
functions $f(x)$ satisfying the following conditions:

\smallskip
1) $f(x)=\sum\limits_{n=0}^{\infty}a_n\cos(2\pi nx)$;

\smallskip
2) $a_n\ge0$ for all $n=0,1,2,\ldots$;

\smallskip
3) $f(0)=\sum\limits_{n=0}^{\infty}a_n=1$;

\smallskip
4) $f(x)=0$ for $h\le|x|\le1/2$.

\smallskip
It is required to evaluate the quantity
\[
A(h)=\sup_{f\in K(h)}a_0=\sup_{f\in K(h)}\int_{-h}^hf(x)\,dx.
\]
\end{Turanpro}

Let $p,\,q\in\N$, $2p\le q$, $(p,q)=1$. In 1972 S.B.~Stechkin (see \cite{GorMan}) solved the problem
for $p=1$, $q=2,3,\ldots$ ($A(1/q)=1/q$). In \cite{GorMan}, the value of $A(h)$ was
calculated for $p=2,\,3$ and $q=2p+1$.

Consider sets
\[
K_{p,q}^0=\{p,p+1,\ldots,q-p\},\qquad
K_{p,q}=q\Z_++K_{p,q}^0.
\]

\begin{thm}\label{t1}
For $p=1$\tu, $q=2,3,\ldots$
\[
\d(K_{1,q}^0)=\d(K_{1,q})=A(1/q)=\frac1q.
\]

For $p=2$\tu, $q=3,5,\ldots$
\[
\d(K_{2,q}^0)=\d(K_{2,q})=A(2/q)=\frac{1+\cos(\pi/q)}{q\cos(\pi/q)}.
\]

For $p=3$, $q=7,8,10,11,\dots$
\[
\d(K_{3,q}^0)=\d(K_{3,q})=A(3/q)=\frac{1}{q}\,\Bigl(1+\frac{1-2(\cos(2\pi
r_0/q)+\cos(2\pi(r_0+1)/q))}{1+2\cos(2\pi r_0/q)\cos(2\pi(r_0+1)/q)}\Bigr),
\]
where $r_0=[q/3]$ is an integer part of $q/3$.

For $q=2p+1$\tu, $p=1,2,\dots$
\[
\d(K_{p,2p+1}^0)=\d(K_{p,2p+1})=A(p/(2p+1))=\frac{\cos(\pi/(2p+1))}{1+\cos(\pi/(2p+1))}.
\]
\end{thm}

Since $K_{1,q}^0=K_q^0$, $K_{1,q}=K_q$, $K_{2,3}^0=\{2,3\}$
and
\[
A(1/q)=\frac1q,\qquad A(2/5)=\frac{1+\cos(\pi/5)}{5\cos(\pi/5)}=\frac{\cos(\pi/5)}{1+\cos(\pi/5)},
\]
then by theorem~\ref{t1}, so that we get examples ~\ref{p1} and~\ref{p2}.

\begin{proof}[\bf Proof of theorem~\ref{t1} in the case $p=3$]
Let $f(x)$ be arbitrary function of class $K(p/q)$, and let $T^{(\e)}\in\T(K_{p,q})$ be
polynomial such that
\begin{equation}\label{f1}
T_0^{(\e)}\le\d(K_{p,q})+\e,
\end{equation}
where $\e>0$ is a little number. Since $T_k^{(\e)}=0$ for $k\in\N\setminus
K_{p,q}$, $f(k/q)=0$ for $k=q\nu+k'$ ($\nu\in\Z_+$, $k'=p,\ldots,q-p$), i.e. for
$k\in K_{p,q}$ (by definition of set $K_{p,q}$), it follows that
\[
T_0^{(\e)}=\sum_{\substack{k\in K_{p,q}\cup\{0\}\\k\le\deg T^{(\e)}}}T_k^{(\e)}f(k/q)=
\sum_{n=0}^\infty a_n\sum_{\substack{k\in K_{p,q}\cup\{0\}\\k\le\deg T^{(\e)}}}
T_k^{(\e)}c(nk/q)=\sum_{n=0}^\infty a_n T^{(\e)}(nk/q).
\]
Hence by nonnegativeness of $a_n$ $\forall n\in\N$ and $T^{(\e)}(x)$ $\forall x\in\R$,
by equality $T^{(\e)}(0)=1$ and inequality \eqref{f1}, so that
\[
a_0\le \d(K_{p,q})+\e.
\]
Since this inequality is established for arbitrary function $f\in K(p/q)$, it
follows that
\[
\sup_{f\in K(p/q)}a_0=A(p/q)\le\d(K_{p,q})+\e.
\]
Hence we obtain lower estimate for value $\d(K_{p,q})$ as $\e\to0$:
\[
A(p/q)\le\d(K_{p,q}).
\]

Let us show that $\d(K_{p,q}^0)\le A(p/q)$ for $p=3$, $q=7,8,10,11,\ldots$ Suppose
\[
\G(\nu)=\g_0+\g_1\cos(2\pi r_0\nu/q)+\g_2\cos(2\pi(r_0+1)\nu/q),\qquad r_0=[q/3],\quad
\nu\in\mathbb Z,
\]
where the coefficients $\g_i$ choose from the equations
\[
\G(0)=1,\ \G(1)=\G(2)=0\iff
\left\{\begin{aligned}
&\g_0+\g_1+\g_2=1,\\
&\g_0+\g_1\cos(2\pi r_0/q)+\g_2\cos(2\pi(r_0+1)/q)=0,\\
&\g_0+\g_1\cos(4\pi r_0/q)+\g_2\cos(4\pi(r_0+1)/q)=0.
\end{aligned}\right.
\]

In \cite{GorMan} prove that $\g_i>0$, $i=0,1,2$, and
\begin{equation}\label{f3}
\frac1{q\g_0}=A(3/q).
\end{equation}

Let $F(x)=F_q(x)$ is Fej\'er's polynomial
\[
F(x)=\sum_{\nu=0}^{q-1}F_\nu\cos(2\pi\nu x)=
\frac1q\biggl(1+2\sum_{\nu=1}^{q-1}\Bigl(1-\frac\nu q\Bigr)
\cos(2\pi\nu x)\biggr)=\Bigl(\frac{\sin(\pi qx)}{q\sin(\pi x)}\Bigr)^2.
\]

Consider polynomial $T^*(x)$
\begin{multline}\label{f2}
T^*(x)=F(x)+\frac{\g_{1}}{2\g_{0}}\,(F(x+r_0/q)+F(x-r_0/q))
\br+
\frac{\g_{2}}{2\g_{0}}\,(F(x+(r_0+1)/q)+F(x-(r_0+1)/q))
\br=
\sum_{\nu=0}^{q-1}\g_0^{-1}F_\nu(\g_0+\g_1\cos(2\pi r_0\nu/q)+\g_2\cos(2\pi(r_0+1)\nu/q))
\cos(2\pi\nu x)
\br=
\sum_{\nu=0}^{q-1}\g_0^{-1}\G(\nu)F_\nu\cos(2\pi\nu x).
\end{multline}

The polynomial $T^*(x)$ satisfies conditions 1)--3) of set $\T(K_{3,q}^0)$.

Since $\G(\nu)=\G(q-\nu)=0$ as $\nu=1,2$, $\G(0)=1$, $F_0=1/q$, it follows that
\[
T^*(x)=\frac1{q\g_0}+\sum_{\nu=3}^{q-3}\g_0^{-1}\G(\nu)F_\nu\cos(2\pi\nu x)=
T_0^*+\sum_{k\in K_{3,q}^0}T_k^*\cos(2\pi\nu x).
\]
Therefore the polynomial $T^*(x)$ satisfies condition 1) of set $\T(K_{3,q}^0)$.

All coefficients $\g_i>0$ and $A(3/q)=(q\g_0)^{-1}$ \eqref{f3}. Hence, by positiveness of
polynomial $F(x)$ and definition of polynomial $T^*(x)$ \eqref{f2}, so that $T^*(x)\ge0$
$\forall x\in\R$. The condition 2) is satisfied.

Now we verify the condition 3). We have
\[
T^*(0)=F(0)+\frac{\g_{1}}{\g_{0}}\,F(r_0/q)+\frac{\g_{2}}{\g_{0}}\,F((r_0+1)/q).
\]
However $F(\nu/q)=0$ for $\nu=1,2,\ldots,q-1$, and points $r_0$, $r_0+1$
are integers from an interval $(0,q/2)$. Therefore $F(r_0/q)=F((r_0+1)/q)=0$ and
\[
T^*(0)=F(0)=1.
\]

Thus $T^*(x)$ belongs to set $\T(K_{3,q}^0)$ and
\[
T_0^*=\frac1{q\g_0}=A(3/q).
\]
Hence upper estimate is $\d(K_{3,q}^0)\le T_0^*=(q\g_0)^{-1}=A(3/q)$.

Finally using $K_{p,q}^0\subset K_{p,q}$ and property 1) of $\d(K)$, we get
estimates
\[
A(p/q)\le\d(K_{p,q})\le\d(K_{p,q}^0)\le A(p/q),
\]
i.e.
\[
\d(K_{p,q}^0)=\d(K_{p,q})=A(p/q).
\]
Polynomial $T^*(x)$ belongs to set $\T(K_{p,q}^0)$ and belongs to set
$\T(K_{p,q})\supset\T(K_{p,q}^0)$. It is extremal polynomial. This completes the
proof of theorem.
\end{proof}

\end{document}